\documentclass[11pt]{article}

\usepackage[T1]{fontenc}
\usepackage{lmodern}
\usepackage[margin=1in]{geometry}
\usepackage{amsmath,amssymb,amsthm,mathtools}
\usepackage{enumitem}
\usepackage{booktabs,array}
\usepackage{xcolor}
\usepackage[protrusion=true,expansion=false]{microtype}
\usepackage{fvextra}
\usepackage{xurl}
\usepackage{authblk}

\usepackage{tikz}
\usetikzlibrary{arrows.meta,fit,positioning}

\usepackage[colorlinks=true,linkcolor=blue,citecolor=blue,urlcolor=blue]{hyperref}
\hypersetup{
  pdftitle={A Lean Formalization of the Hamilton--Perelman Proof of Thurston's Geometrization Conjecture},
  pdfauthor={Yuan Liao, Ziyang Qin, Jui-Hui Chung, Ayush Khaitan, Zhengyang Zhang, 
  Hongzhou Lin, 
  Zijian Zhang, Ligeng Zhu, Peng Lu, and Bennett Chow}
}

\theoremstyle{plain}
\newtheorem{theorem}{Theorem}

\theoremstyle{definition}

\theoremstyle{remark}

\newcommand{\Ric}{\operatorname{Ric}}

\newcommand{\Vol}{\operatorname{Vol}}

\newcommand{\releasetag}{v0.1.5}
\newcommand{\releasecommit}{bafefb7ce5ce7e40417641d67e6e95727b7bba1f}
\newcommand{\releaseblobbase}{https://github.com/qinz1yang/differential-geometry/blob/\releasecommit/DifferentialGeometry/}
\newcommand{\releasetreebase}{https://github.com/qinz1yang/differential-geometry/tree/\releasecommit/DifferentialGeometry/}
\newcommand{\declareleanfile}[3]{%
  \expandafter\def\csname leanfile@#1\endcsname{#2}%
  \expandafter\def\csname leanline@#1\endcsname{#3}}
\newcommand{\leanref}[1]{%
  \ifcsname leanfile@#1\endcsname
    \href{\releaseblobbase\csname leanfile@#1\endcsname\#L\csname leanline@#1\endcsname}{\nolinkurl{#1}}%
  \else
    \nolinkurl{#1}%
  \fi}
\newcommand{\leanrefas}[2]{%
  \ifcsname leanfile@#1\endcsname
    \href{\releaseblobbase\csname leanfile@#1\endcsname\#L\csname leanline@#1\endcsname}{#2}%
  \else
    #2%
  \fi}
\newcommand{\leanfile}[1]{\nolinkurl{#1}}

\newcommand{\rfdir}[1]{\href{\releasetreebase#1}{\leanfile{#1/}}}
\newcommand{\releasebranch}{paper/geometrization-v0.1.5}
\newcommand{\Nil}{\mathrm{Nil}}
\newcommand{\Sol}{\mathrm{Sol}}
\declareleanfile{geometrization}{Geometry/Flow/RicciFlow/LongTime/Geometrization.lean}{49}
\declareleanfile{geometrization_certificate}{Geometry/Flow/RicciFlow/LongTime/Geometrization.lean}{33}
\declareleanfile{smooth_geometrization_conjecture}{Geometry/Flow/RicciFlow/LongTime/Geometrization.lean}{112}
\declareleanfile{ConnectedClosedOrientedManifold}{Topology/Manifold/ClosedOriented.lean}{68}
\declareleanfile{IsPrime}{Topology/ThreeManifold/Geometrization/Prime.lean}{11}
\declareleanfile{TorusDecomposition}{Topology/ThreeManifold/TorusCut/Decomposition.lean}{18}
\declareleanfile{exists_surgery_with_decaying_accuracy_enhanced}{Geometry/Flow/RicciFlow/LongTime/Ch11/A12Enhanced.lean}{36}
\declareleanfile{late_sequence_tests}{Geometry/Flow/RicciFlow/LongTime/Ch11/RouteWLateSequenceCoreC11R.lean}{39}
\declareleanfile{components_geometrize_of_late_sequence_tests}{Geometry/Flow/RicciFlow/LongTime/LateDecomposition.lean}{136}
\declareleanfile{incompressible_of_shifted_area_barriers}{Geometry/Flow/RicciFlow/LongTime/CuspIncompressibility.lean}{13}
\declareleanfile{not_nonnegative_area_upper_barriers_pi_shift}{Analysis/ODE/AreaUpperBarrier.lean}{102}
\declareleanfile{geometrizes_of_late_slice_supply}{Geometry/Flow/RicciFlow/LongTime/RegularSlice.lean}{91}
\declareleanfile{geometrizes_of_raw_empty}{Geometry/Flow/RicciFlow/Surgery/Geometrization/ClassifiedHistory.lean}{24}
\declareleanfile{exists_primeDecomposition}{Topology/ThreeManifold/PrimeDecomposition/Existence.lean}{43}
\declareleanfile{primeSeam}{Topology/ThreeManifold/Geometrization/SmoothTorusReconstruction.lean}{90}
\declareleanfile{HasThurstonAtlas}{Geometry/Thurston/Atlas.lean}{21}
\makeatletter
\g@addto@macro{\UrlBreaks}{\do\_\do\.\do\/}
\makeatother
\title{A Lean Formalization of the Hamilton--Perelman Proof of Thurston's Geometrization Conjecture\thanks{Yuan Liao and Ziyang Qin contributed equally.}}
\author[1]{Yuan Liao}
\author[2]{Ziyang Qin}
\author[3]{Jui-Hui Chung}
\author[3]{Ayush Khaitan}
\author[4]{Zhengyang Zhang}
\author[5]{Hongzhou Lin}
\author[4,6]{Zijian Zhang}
\author[6]{Ligeng Zhu}
\author[7]{Peng Lu}
\author[1]{Bennett Chow}

\affil[1]{Department of Mathematics, University of California San Diego}
\affil[2]{Department of Mathematics, Cornell University}
\affil[3]{Princeton AI, Princeton University}
\affil[4]{Tsinghua University}
\affil[5]{Amazon AGI}
\affil[6]{NVIDIA Research}
\affil[7]{Department of Mathematics, University of Oregon}
\date{}

\begin{document}
\maketitle

We formalize the Hamilton--Perelman Ricci-flow proof of
Thurston's geometrization conjecture. The public Lean declaration is reproduced below;
its geometric interpretation follows, and verification is recorded in
Section~\ref{sec:verification}.

\section{Main theorems}

The endpoint is
\leanrefas{geometrization}{\texttt{GC.Endpoint.geometrization}}.
Its declaration type is reproduced here with the proof omitted:
\begin{Verbatim}[fontsize=\footnotesize,breaklines=true,commandchars=\\\{\}]
theorem GC.Endpoint.geometrization (M : ConnectedClosedOrientedManifold.\{u\} 3) :
    \ensuremath{\exists} factors : List (ConnectedClosedOrientedManifold.\{u\} 3),
      factors \ensuremath{\ne} [] \ensuremath{\land}
      (\ensuremath{\forall} P \ensuremath{\in} factors, IsPrime P) \ensuremath{\land}
      Nonempty (ClosedOrientedManifold.OrientedDiffeomorph
        (finiteConnectedSum factors).toClosedOrientedManifold M.toClosedOrientedManifold) \ensuremath{\land}
      \ensuremath{\forall} P \ensuremath{\in} factors, \ensuremath{\exists} D : GC.Topology.TorusDecomposition P,
        let T := D.reconstructionAtlas.torusInPrime D.reconstruction
        (\ensuremath{\forall} i, ContMDiff torusModel (\ensuremath{\mathcal R} 3) \ensuremath{\infty} (T i) \ensuremath{\land} _root_.Topology.IsEmbedding (T i)) \ensuremath{\land}
        Pairwise (fun i j => Disjoint (Set.range (T i)) (Set.range (T j))) \ensuremath{\land}
        (\ensuremath{\forall} i x, Function.Injective (FundamentalGroup.map (T i) x)) \ensuremath{\land}
        \ensuremath{\forall} j : Fin D.components.count,
          let := Manifold.interiorChartedSpace D.carrier.model \ensuremath{\infty}
            (M := D.carrier.pieceInterior (D.components.piece j))
          let := Manifold.interiorIsManifold D.carrier.model \ensuremath{\infty}
            (M := D.carrier.pieceInterior (D.components.piece j))
          \ensuremath{\exists} (k : GC.Geometry.ThurstonModel)
            (g : SmoothRiemannianMetric (\ensuremath{\mathcal R} 3)
              (D.carrier.pieceInterior (D.components.piece j))),
            RiemannianMetricComplete g \ensuremath{\land} GC.Geometry.HasThurstonAtlas g k \ensuremath{\land}
              (k = .hyperbolic \ensuremath{\to} Integral.Measure.riemannianVolumeMeasure (\ensuremath{\mathcal R} 3)
                (D.carrier.pieceInterior (D.components.piece j)) g Set.univ < \ensuremath{\top})
\end{Verbatim}

\clearpage
The geometric content of this declaration is the following.

\begin{theorem}[Geometrization; Hamilton--Perelman]
\label{thm:geometrization}
Let \(M\) be a connected, closed, oriented smooth three-manifold.
There is an orientation-preserving diffeomorphism
\[
 M\cong M_1\mathbin{\#}\cdots\mathbin{\#}M_k,\qquad k\geq1,
\]
where each \(M_i\) is a connected, closed, oriented smooth prime
three-manifold. In each \(M_i\) there is a finite, possibly empty family of
pairwise disjoint, smoothly embedded, two-sided incompressible tori.
For every piece \(N\) obtained by cutting along these tori, its interior
\(N^\circ\), with the inherited smooth structure, admits a smooth
Riemannian metric \(g\) and a model \(X\) among the eight Thurston geometries
\[
 \mathbb S^3,\quad \mathbb E^3,\quad \mathbb H^3,\quad
 \mathbb S^2\!\times\!\mathbb R,\quad \mathbb H^2\!\times\!\mathbb R,\quad
 \widetilde{\mathrm{SL}_2(\mathbb R)},\quad \Nil,\quad \Sol,
\]
such that \(g\) is complete and has an atlas of local isometries to \(X\).
If \(X=\mathbb H^3\), then \(\Vol_g(N^\circ)<\infty\).
\end{theorem}

The underlying topological framework has two stages:
Kneser's existence theorem~\cite{Kneser1929} and Milnor's
uniqueness theorem~\cite{Milnor1962} describe the decomposition
into prime connected summands, while the
Jaco--Shalen--Johannson theory~\cite{JacoShalen1979,Johannson1979}
decomposes the irreducible summands along incompressible tori
into Seifert-fibered and atoroidal pieces.
Theorem~\ref{thm:geometrization} asserts the existence of a prime
connected-sum decomposition and, in each factor, a decomposition
along incompressible tori such that the interior of every resulting
piece admits a complete Thurston geometry; it does not require the cutting
tori to form the canonical JSJ family.

Here closed means compact without boundary. Incompressibility means that
inclusion of each torus induces an injection on fundamental groups.
Our \leanrefas{IsPrime}{prime convention} allows \(\mathbb S^3\) as a factor,
so that the nonempty connected sum also includes the sphere itself.
We do not require the torus family to be minimal and make no uniqueness
assertion. We cut only along tori and impose finite volume only in the
hyperbolic case; nonhyperbolic interiors may have infinite volume.
For the torus-only convention and the role of the twisted real-line bundle
over the Klein bottle, see \cite[Section~2 and Remark~2.2]{Lott2010}.

The public Lean endpoint is
\leanrefas{geometrization}{\texttt{GC.Endpoint.geometrization}}.
Its \leanrefas{ConnectedClosedOrientedManifold}{input type} bundles a compact
Hausdorff smooth three-manifold, connectedness, and an orientation.
No metric, surgery history, irreducibility assumption, or geometric
classification is required from the caller. The metric used to begin the
proof is obtained from compactness.

The conclusion returns a finite nonempty list of prime factors and an
oriented diffeomorphism from their connected sum to the original manifold.
For each factor it returns a \leanrefas{TorusDecomposition}{torus decomposition}, including the
cut carrier, its connected pieces, the paired boundary tori, and a smooth
reconstruction of that same factor. The returned tori are explicitly
smooth, embedded, pairwise disjoint, and injective on fundamental groups;
the reconstruction also supplies their
\leanrefas{primeSeam}{smooth two-sided collars}.
For each actual piece interior, the conclusion then exposes one model
and one smooth metric, together with completeness, a
\leanrefas{HasThurstonAtlas}{local isometric atlas},
and the hyperbolic finite-volume implication for that same metric.
The model is uniform on the piece.

The declaration \leanref{geometrization_certificate} retains the
construction data. The unbundled statement
\leanref{smooth_geometrization_conjecture} expresses the same certificate
conclusion for an ordinary smooth three-manifold equipped with an orientation.
These are formulations of the same result, rather than additional
classification hypotheses.

\section{Mathematical route}

In this section we summarize the mathematical route followed by our Lean formalization.
The proof 
of Thurston's geometrization conjecture \cite{Thurston1982}
follows Hamilton's Ricci-flow program
\cite{Hamilton1982,Hamilton1993Harnack,Hamilton1995Singularities,Hamilton1997,Hamilton1999} and Perelman's completion of it \cite{Perelman2002,Perelman2003Surgery}.
Morgan--Tian \cite{MorganTian2007,MorganTian2014} and
Kleiner--Lott \cite{KleinerLott2008} both provide global organizations and detailed accounts of
Perelman's arguments. The Perelman collapse step \cite[\S 7.4]{Perelman2003Surgery}
(generalizing Cheeger--Gromov collapsing theory \cite{CheegerGromov1986})
uses the local-collapse reworking
of Kleiner--Lott \cite{KleinerLott2014} (cf. Perelman \cite{Perelman1991} and~Shioya--Yamaguchi \cite{ShioyaYamaguchi2005}), including its boundary version
and application to geometrization in Sections~16--17. The passage from
graph manifolds to geometric pieces uses Seifert-fibered geometry \cite{JacoShalen1979} and
the eight model geometries \cite{Scott1983,Thurston1997}.

Start with a smooth Riemannian metric \(g_0\) on \(M\).
Between surgery times the metric satisfies
\[
 \partial_t g=-2\Ric_g.
\]
The short-time existence, curvature estimates, noncollapsing, canonical neighborhoods,
and surgery constructions are shared with the Poincar\'e conjecture development~\cite{ChowLiaoQin2026,QinLiaoKhaitanChow2026} (cf.~\cite{ZhangEtAl2026}).
For geometrization, the flow must supply information at arbitrarily late
times whenever components survive. The \leanrefas{exists_surgery_with_decaying_accuracy_enhanced}{long-time surgery existence theorem}
produces a
surgery flow whose accuracy function \(\delta(t)\) decreases to zero,
together with the analytic estimates used in the long-time argument.
The construction keeps finite surgery histories on finite observation
intervals, with identifications of the initial manifold and the surviving
and discarded components.

At a regular time \(t>0\), one studies the rescaled metrics
\(t^{-1}g(t)\); the normalization \((4t)^{-1}g(t)\) is used for
hyperbolic curvature \(-1\).
The thick-region analysis gives compact cores of complete finite-volume
hyperbolic manifolds, with cusp tori separating them from the remaining
region. The relevant cuts and their smooth reconstruction are carried
along with the estimates. They are used both by the subsequent collapse
argument and by the proof that the cutting tori are incompressible.

On the remaining region, local volume collapse at the curvature scale is
combined with curvature-derivative control. For pieces with boundary,
the hypotheses include nearly cuspidal boundary collars. The collapse
theorem yields graph-manifold pieces, with closed geometric alternatives
handled separately. In particular, the formal late-piece classification
also includes a closed nonnegative-curvature case; it is not necessary
to recognize every such component as a graph manifold.
Refining the graph-manifold presentation into geometric pieces, while
retaining its boundary identifications, supplies the required local
Thurston geometries. The resulting geometric metrics are complete on the interiors of
the refined pieces.

The analytic inputs are organized as
\leanrefas{late_sequence_tests}{tests on arbitrarily late sequences of regular slices}.
For every sequence of nonempty regular slices with times tending to
infinity, the tests eventually provide torus cuts and an exterior-area
obstruction to compressibility. Each piece satisfies the hyperbolic,
collapsed, or closed nonnegative-curvature alternative. The derivative
control determines a collapse threshold. If non-geometrizable components
occurred at arbitrarily late times, one could choose a sequence of such
slices; the threshold and the eventual tests would geometrize a member
of this sequence, a contradiction. This is formalized by the
\leanrefas{components_geometrize_of_late_sequence_tests}{late-sequence geometrization theorem}.

For incompressibility, suppose that one of the relevant torus inclusions
fails to induce an injection on fundamental groups. 
Hamilton's incompressibility argument~\cite{Hamilton1999},
drawing on the minimal-surface methods of Meeks--Yau
\cite{MeeksYau1982} and Schoen--Yau~\cite{SchoenYau1979}
and Mostow--Prasad rigidity \cite{Mostow1968,Prasad1973},
is implemented through a continuous nonnegative area function 
\(A:[T,\infty)\to\mathbb R\), for some \(T\geq0\), and a constant
\(c>0\), satisfying the upper-barrier inequality
\[
 A'(t)\leq\frac{3}{4(t+c)}A(t)-\pi.
\]
This notation is understood in the local smooth upper-barrier sense;
ordinary differentiability of \(A\) is not assumed. The moving boundary
and surgery comparisons are part of the construction of this barrier.
Using \(A\geq0\), the formal comparison shows that
\[
 F(t):=\frac{A(t)}{t+c}+\pi\log(t+c)
\]
is nonincreasing. Consequently, for \(t\geq T\),
\[
 0\leq\frac{A(t)}{t+c}
 \leq\frac{A(T)}{T+c}-\pi\log\!\left(\frac{t+c}{T+c}\right).
\]
The right-hand side becomes negative for large \(t\).
The contradiction is formalized in the
\leanrefas{not_nonnegative_area_upper_barriers_pi_shift}{upper-barrier area lemma}; applying it to
the exterior areas proves
\leanrefas{incompressible_of_shifted_area_barriers}{incompressibility}.
Thus the tori used in the geometric decomposition are incompressible in
the reconstructed manifold, rather than merely in an auxiliary model.

Once the sufficiently late nonempty components geometrize, the finite
surgery histories reconstruct the original manifold. The reconstruction
includes the standard factors discarded during surgery and the connected
sums created when cuts are reversed. It also covers an empty terminal
flow, through the \leanrefas{geometrizes_of_raw_empty}{empty-terminal reconstruction theorem}.
Consequently the argument does not assume that every initial manifold
becomes extinct in finite time. The \leanrefas{geometrizes_of_late_slice_supply}{late-slice reconstruction theorem}
combines the late-component
supply with this reconstruction.
\leanrefas{exists_primeDecomposition}{Prime decomposition}
and geometric reconstruction produce a certificate for the original
\(M\); the public endpoint then exposes the factors, actual torus cuts,
and complete metrics stated in Theorem~\ref{thm:geometrization}.
A fuller account of the analytic and topological constructions will
appear elsewhere.

\section{Agent-assisted formalization with Humanize}

We use the Humanize RLAR (Ralph Loop with Agentic Review; see \href{https://humanfia.ai/flows/rlar}{https://humanfia.ai/flows/rlar}) workflow~\cite{humanize2026} to automate Lean formalization while retaining human control. Humanize separates implementation from review. A builder agent develops the code in successive rounds, while an independent reviewer checks progress against the plan and inspects the changes for errors. This structure supports long-running formalization by recording evidence of progress and unresolved proof obligations, alongside compiler diagnostics and reviewer feedback.

The formal-proof workflow combines kernel checking, statement-preservation checks, and independent proof verification. These checks complement, rather than replace, mathematical review of the definitions and hypotheses. Human authors \textit{remain responsible} for the intended theorem statements, the mathematical interpretation of the formalization, and the final verification claims.

We use GPT-6-Astra as the reviewer. For the builder, we choose among GPT-6-Luna, GPT-6-Sol, GPT-6-Astra, Claude-Opus-5.5, and DeepSeek-V4.1-Flash based on their observed generation speed in tokens per second (TPS).

\section{Formal artifact and verification}
\label{sec:verification}

The formalization uses Lean~4 \cite{Lean2021} and Mathlib
\cite{Mathlib2020}. We fix the source revision below, also released as
\href{https://github.com/qinz1yang/differential-geometry/releases/tag/\releasetag}{\releasetag\ (preview)}.
The paper branch records this revision; all source links use its complete
commit SHA.

\begin{center}
\small
\begin{tabular}{@{}ll@{}}
\toprule
Repository & \url{https://github.com/qinz1yang/differential-geometry}\\
Paper branch & \href{https://github.com/qinz1yang/differential-geometry/tree/\releasebranch}{\texttt{\releasebranch}}\\
Tag & \texttt{\releasetag}\ (preview)\\
Commit & \href{https://github.com/qinz1yang/differential-geometry/commit/\releasecommit}{\texttt{\releasecommit}}\\
Lean / Mathlib & \texttt{v4.35.0-rc3} / \texttt{v4.35.0-rc3}\\
\bottomrule
\end{tabular}
\end{center}

The toolchain and package revisions are pinned in
\texttt{lean-toolchain} and \texttt{lake-manifest.json}.
Vendored source and adaptation records are fixed by the same commit.
The project retains the source attribution and licenses of the
Poincar\'e conjecture development, including selected contributions from Armstrong and Kempe's
\texttt{DeGiorgi} library for De Giorgi--Nash--Moser theory,
Jordan--Sch\"onflies, Classification of Compact Surfaces, Tau Ceti,
and the canonical-topology development; the preserved notices and
source headers identify the respective adaptations. These uses concern
selected developments, not the complete theorem collections of those
projects. The root \href{https://github.com/qinz1yang/differential-geometry/blob/\releasecommit/NOTICE}{\texttt{NOTICE}}
and the \rfdir{External} directory provide the corresponding attribution
and modification records.

From the repository root at the pinned revision, build the endpoint and
the complete library with:
\begin{Verbatim}[fontsize=\footnotesize,breaklines=true]
lake build DifferentialGeometry.Geometry.Flow.RicciFlow.LongTime.Geometrization
lake build DifferentialGeometry
\end{Verbatim}

The recorded verification includes successful endpoint and complete-library
builds, a separate consumer of the expanded public statement, and
transitive axiom checks. Declaration-linter checks were applied to the
four public geometrization declarations listed below.
The source and dependency contents used for those checks are preserved
through the history cleanup and release integration; this transfers the
recorded evidence by source identity.

To inspect the public statement and print its transitive axiom closure,
put the following in a separate file, for example
\texttt{/tmp/GeometrizationAxioms.lean}:
\begin{Verbatim}[fontsize=\footnotesize,breaklines=true]
import DifferentialGeometry.Geometry.Flow.RicciFlow.LongTime.Geometrization

#check GC.Endpoint.geometrization
#print axioms GC.Endpoint.geometrization
#print axioms GC.Endpoint.geometrization_certificate
#print axioms GC.Endpoint.geometrization_conjecture
#print axioms GC.Endpoint.smooth_geometrization_conjecture
\end{Verbatim}
Run it from the repository root with:
\begin{Verbatim}[fontsize=\footnotesize,breaklines=true]
lake env lean --root=/tmp /tmp/GeometrizationAxioms.lean
\end{Verbatim}
The four axiom queries yield the same set:
\[
 [\texttt{propext},\ \texttt{Classical.choice},\ \texttt{Quot.sound}].
\]
These proof closures contain neither \texttt{sorryAx} nor a
project-specific axiom. These checks concern the named declarations and
their dependencies; their mathematical interpretation also requires
review of the definitions and hypotheses.

\section{Acknowledgements}
A.K.~acknowledges funding from the DARPA expMath Program, as well as the AMS-Simons Travel Grant.
B.C.~thanks Jingwei Xia for her constant support and OpenAI for providing access to OpenAI models through the ChatGPT for Academic Researchers program.

\end{document}